\documentclass[12pt,psamsfonts]{amsart}

\usepackage{amsmath,amssymb}
\usepackage{graphicx}
\usepackage{color}
\usepackage{psfrag} % See the section on figures.
\usepackage{overpic} % See the section on figures.
\usepackage{array} % Function notation
\usepackage{subfigure}
 
\newtheorem{theorem}{Theorem}
\newtheorem{proposition}[theorem]{Proposition}

\newtheorem{lemma}[theorem]{Lemma}
\newtheorem {remark}[theorem]{Remark}
\newtheorem{definition}[theorem]{Definition}

\title[Analytic nilpotent centers on center manifolds]{Analytic nilpotent centers on center manifolds}
\author[L. Queiroz and C. Pessoa]{}


\begin{document}
 \maketitle

\centerline{\scshape  Claudio Pessoa,  \; Lucas Queiroz}
\medskip

{\footnotesize \centerline{Universidade Estadual Paulista (UNESP), Instituto de Bioci\^encias Letras e Ci\^encias Exatas,} \centerline{R. Cristov\~ao Colombo, 2265, 15.054-000, S. J. Rio Preto, SP, Brasil }
\centerline{\email{c.pessoa@unesp.br} and \email{lucas.queiroz@unesp.br}}}

\medskip

\bigskip

\begin{quote}{\normalfont\fontsize{8}{10}\selectfont
{\bfseries Abstract.}
Consider analytical three-dimensional differential systems having a singular point at the origin such that its linear part is $y\partial_x-\lambda z\partial_z$ for some $\lambda\neq 0$. The restriction of such systems to a Center Manifold has a nilpotent singular point at the origin. We prove that if the restricted system has an analytic nilpotent center at the origin, with Andreev number $2$, then the three-dimensional system admits a formal inverse Jacobi multiplier. We also prove that nilpotent centers of three-dimensional systems, on analytic center manifolds, are limits of Hopf-type centers. We use these results to solve the center problem for some three-dimensional systems without restricting the system to a parametrization of the center manifold.
\par}
\end{quote}

\section{Introduction}

Consider the analytical vector fields $X$ in $\mathbb{R}^3$ having a singular point $p$ such that $DX(p)$ has two zero eigenvalues and one real eigenvalue $\lambda\neq 0$ and the rank of $DX(p)$ is $2$. The differential system associated to these vector fields can always be written, by an affine change of variables, in the form:
\begin{equation}\label{eq1}
	\begin{array}{lr}
		\dot{x}=y+P(x,y,z),\\
		\dot{y}=\hspace{0.7cm}Q(x,y,z),\\
		\dot{z}=-\lambda z + R(x,y,z),
	\end{array}
\end{equation}
where $P,Q,R$ are analytic functions and $j^1P(0)=j^1Q(0)=j^1R(0)=0$. By the Center Manifold Theorem, we know that for every $r\geqslant 1$, there exists a bidimensional $C^r$-manifold invariant by system \eqref{eq1} tangent to the $xy$-plane at the origin (see \cite{Kelley,Sijbrand}).

%having a singular point $p$ such that $DX(p)$ has the following Jordan canonical form
%$$\begin{small}
%	\left(\begin{array}{lcr}
%		0 & 1 & 0 \\
%		0 & 0 & 0 \\
%		0 & 0 & -\lambda
%	\end{array}
%	\right),
%\end{small}$$
%where $\lambda\neq 0$. We can assume that the singular point is located at the origin. The associated differential system can be written, by a linear change of variables, in the following form:
%\begin{equation}\label{eq1}
%	\begin{array}{lr}
%		\dot{x}=y+P(x,y,z),\\
%		\dot{y}=\hspace{0.7cm}Q(x,y,z),\\
%		\dot{z}=-\lambda z + R(x,y,z),
%	\end{array}
%\end{equation}
%where $P,Q,R$ are analytic functions and $j^1P(0)=j^1Q(0)=j^1R(0)=0$. For such systems, the following theorem holds.

%\begin{theorem}[Center Manifold Theorem]\label{TeoCenterManifold} Consider system \eqref{eq1}. There exists an invariant bidimensional $C^r$-manifold tangent to the $xy$-plane at the origin for every $r\geqslant 1$.
%\end{theorem}

These manifolds are called \emph{Center Manifolds} and neither the analiticity nor the uniqueness of the center manifold is guaranteed by the Center Manifold Theorem. The restriction of system \eqref{eq1} to a center manifold is a bidimensional system. Since every center manifold is tangent to the $xy$-plane at the origin, the restricted system can be written in the form
\begin{equation}\label{eqPlanar}
	\begin{array}{lr}
		\dot{x}=y+P(x,y),\\
		\dot{y}=Q(x,y),
	\end{array}
\end{equation}
where $j^1P(0)=j^1Q(0)=0$. Thus, the restricted system has a nilpotent singular point at the origin. Therefore, we say that a three-dimensional analytical vector field has a \emph{nilpotent singular point} if its associated differential system can be put in the form \eqref{eq1}.

In the plane, nilpotent singular points are widely studied in the literature \cite{AlvarezGasull1,Llibre,GarciaGineInte}. One of the most important problems in the study of these types of singularities is the so-called \emph{Nilpotent Center Problem} \cite{AlvarezGasull1,Zoladek} which consists of distinguishing whether the origin of \eqref{eqPlanar} is a center or not. In order to study the nilpotent center problem, one must first identify which systems \eqref{eqPlanar} have monodromic singular points at the origin. The following result (see \cite[Theorem 3.5, p. 116]{Llibre} and \cite[Theorem 66, p. 357]{Andronov}) provides a monodromy criterion.

\begin{theorem}[Andreev's Theorem]\label{TeoAndreev}
Let $X$ be the vector field associated to analytic system \eqref{eqPlanar} and the origin be an isolated singular point. Let $y=F(x)$ be the solution of the equation $y+P(x,y)=0$ in a neighborhood of $(0,0)$ and consider $f(x)=Q(x,F(x))$ and $\Phi(x)=\mbox{\rm div}X\vert_{(x,F(x))}$. Then, we can write
	$$f(x)=ax^{\alpha}+O(x^{\alpha+1}),$$
	$$\Phi(x)=bx^{\beta}+O(x^{\beta+1}).$$
	The origin is monodromic if and only if $a<0,\alpha=2n-1$ and one of the following conditions holds:
	\begin{itemize}
		\item[i)]$\beta>n-1$ or $\Phi\equiv 0$;
		\item[ii)]$\beta=n-1$ and $b^2+4an<0$.
	\end{itemize}
\end{theorem}
The positive integer $n$ in the above theorem is called \emph{Andreev number}.
\medskip

Our goal in this work is to study the Nilpotent Center Problem for three-dimensional system \eqref{eq1}, i.e., to determine conditions for the origin to be a nilpotent center on a center manifold. As far as we know, this is a barely explored problem in the literature. The only two works on this subject are the papers \cite{SongWangFeng, Wang2} in which the authors study particular families of system \eqref{eq1}. Their results were obtained by restricting the system to a polynomial approximation of the center manifold. We seek to study the nilpotent center problem without going through the restriction of the system to a center manifold. 

\begin{definition}
	Let $X$ be a vector field defined on an open set $U\subset\mathbb{R}^k$ and $V:U\subset\mathbb{R}^k\to\mathbb{R}$ be a $C^1$-function not locally null such that
	$$XV-V\mbox{\rm div}X\equiv 0.$$
	We say that $V$ is an \emph{inverse Jacobi multiplier} of $X$.
\end{definition}

\begin{remark}
	When $k=2$ in the above definition, the function $V$ is usually called an \emph{inverse integrating factor} of $X$.
\end{remark}

The main results of this paper are the following:

\begin{theorem}\label{TeoJacobin=2}
	Consider system \eqref{eq1} having a monodromic singular point with Andreev number $2$ at the origin. If the origin is an analytic nilpotent center on a center manifold, then there exists a formal inverse Jacobi multiplier $V(x,y,z)$ for system \eqref{eq1} such that $j^{m+1}V(0)=y^mz$ for $m\geqslant 0$.
\end{theorem}

The converse of this theorem is also true when we consider $V(x,y,z)$ satisfying some stronger conditions, for instance, when the center manifold is contained in the zero set $V^{-1}(0)$ and $V(x,y,z)$ is a $C^{\infty}$-function.

\begin{theorem}\label{TeoLlibreLimiteR3}
	Suppose the origin of the system \eqref{eq1} is an analytical nilpotent center on a center manifold. Then there are two functions $F_1(x,y)$, $F_2(x,y)$ analytical at the origin with $j^1F_1(0)=j^1F_2(0)=0$ such that the 1-parameter family
	\begin{equation}\label{eq1Llimit}
		\begin{array}{lr}
			\dot{x}=y+P(x,y,z)+\varepsilon F_1(x,y),\\
			\dot{y}=-\varepsilon x+Q(x,y,z)+\varepsilon F_2(x,y),\\
			\dot{z}=-\lambda z + R(x,y,z).
		\end{array}
	\end{equation}
	has a non-degenerate center at the origin for any $\varepsilon>0$. Also there is an analytic function $f(x,y)$ at the origin with $j^1f(0)=0$ such that $(x-F_1(x,y))\dfrac{\partial f}{\partial y}=F_2(x,y)(1+\dfrac{\partial f}{\partial x})$.
\end{theorem}

This paper is structured in the following way: In Section 2 we recall some preliminary results necessary to develop the study of three-dimensional system \eqref{eq1}. We address the monodromy problem for the restriction of system \eqref{eq1} to a center manifold, introduce a useful Linear Operator which show up frequently in the study of nilpotent singular points in $\mathbb{R}^3$ and  present the formal normal form for system \eqref{eq1}.

In Section 3, we relate the existence of an Inverse Jacobi Multiplier and the Nilpotent Center Problem for system \eqref{eq1}. We exhibit some properties of Inverse Jacobi multipliers and prove that for the vector field $X$ associated to \eqref{eq1} there exists a formal series $V(x,y,z)$ such that $XV-V\mbox{div}X=\sum_{n\geqslant 1}\Lambda_nx^{n-1}z$. In this section we prove Theorem \ref{TeoJacobin=2} and conclude that we can use the quantities $\Lambda_n$ to detect analytic nilpotent centers on center manifolds of system \eqref{eq1}.

In Section 4, we prove Theorem \ref{TeoLlibreLimiteR3} which is an extension of the main result of \cite{LlibreLimite} to three-dimensional systems, i.e. if the origin of \eqref{eq1} is an analytic nilpotent center on a center manifold, than it is a limit of non-degenerate centers of Hopf type on the respective center manifolds. 

Finally, in Section 5, we apply the results of the previous sections to some three-dimensional systems, solving the analytic nilpotent center problem. 

\section{Preliminary results}%\label{SecPreliminary}
In this section, we present some preliminary results on the nilpotent singular points of the three-dimensional system \eqref{eq1}. These results are proven in \cite{QueirozPessoaNillR3}.

\subsection{Monodromy}
We first deal with the lack of analyticity of the restriction of system \eqref{eq1} to a center manifold in the study of the monodromy. In \cite{AlvarezGasull1} the authors have proved Theorem \ref{TeoAndreev}. Following basically the same steps described in their proof it is possible to prove a less restrictive version of Theorem \ref{TeoAndreev}.
%The lack of analyticity of the restricted system \eqref{eq1restricted} apparently poses an obstacle for us to use Theorem \ref{TeoAndreev} to study the monodromy of the singularity on a center manifold (see \cite[Theorem 3.5, p. 116]{Llibre} and \cite[Theorem 66, p. 357]{Andronov}). However, the proof presented in \cite{AlvarezGasull1} uses the concept of the Poincaré index which is not dependent on the analyticity of the vector field. The other important step of the proof uses the \emph{generalized polar coordinates}, first introduced by Lyapunov \cite{LyapunovP} (see also \cite{GasullTorre}), which also do not depend on the analyticity. Thus, by supposing $j^rf(0)\neq 0$, the same arguments can be made and we have a less restrictive version of Theorem \ref{TeoAndreev}.

\begin{theorem}[$C^r$-Andreev's Theorem]\label{TeoAndreevRestrito}
	Let $X$ be the vector field associated to the $C^r$-system, $r\geqslant 3$, given by
	\begin{equation*}%\label{eqPlanarCr}
		\begin{array}{lr}
			\dot{x}=y+X_2(x,y),\\
			\dot{y}=Y_2(x,y),
		\end{array}
	\end{equation*} 	
	where $X_2,Y_2\in C^r$, $j^1X_2(0)=j^1Y_2(0)=0$ and such that the origin is an isolated singular point. Let $y=F(x)$ be the solution of the equation $y+X_2(x,y)=0$ through $(0,0)$  and consider $f(x)=Y_2(x,F(x))$ and $\Phi(x)=\mbox{\rm div}X\vert_{(x,F(x))}$. We can write
	$$f(x)=ax^{\alpha}+O(x^{\alpha+1}),$$
	$$\Phi(x)=bx^{\beta}+O(x^{\beta+1}).$$
	for $\alpha<r$. Suppose that $a\neq 0$, then the origin is monodromic if and only if $a<0, \alpha=2n-1$ and one of the following conditions holds:
	\begin{itemize}
		\item[i)]$\beta>n-1$ or $j^r\Phi(0)\equiv 0$;
		\item[ii)]$\beta=n-1$ and $b^2+4an<0$;
	\end{itemize}
\end{theorem}

%The above theorem provides a tool to study the monodromy of the origin for system \eqref{eq1} on a center manifold. In order to perform this study, we apply Theorem \ref{TeoAndreevRestrito} to the restricted system \eqref{eq1restricted} and check if the function $f$ is not flat, since for any $r\geqslant 3$, there exists a center manifold such that \eqref{eq1restricted} is a $C^r$-system.
%The positive integer $n$ in the statements of Theorems \ref{TeoAndreev} and \ref{TeoAndreevRestrito} plays an important role in the study of nilpotent monodromic singular points. So we define the \emph{Andreev number} of a nilpotent singular point by the number $n$ in function $f(x)=ax^{2n-1}+O(x^{2n})$ described in Theorems \ref{TeoAndreev} and \ref{TeoAndreevRestrito}.

We define the \emph{Andreev number} of a nilpotent singular point by the positive integer $n$ in function $f(x)=ax^{2n-1}+O(x^{2n})$ from Theorems \ref{TeoAndreev} and \ref{TeoAndreevRestrito}. It was proven in \cite{GarciaFiI} that the Andreev number is invariant by analytical and formal orbital equivalence, i.e. via analytical and formal diffeomorphisms and time rescalings. Furthermore, the monodromy conditions (i) and (ii) in Theorem \ref{TeoAndreevRestrito} are invariant by local diffeomorphisms.

\begin{proposition}\label{Propobeta=n-1invariante}
	For $C^r$ system \eqref{eqPlanar} having an isolated monodromic nilpotent singular point at the origin, with Andreev number $n$ such that $2n-1<r$, the monodromy conditions (i) and (ii) in Theorem \ref{TeoAndreevRestrito} are invariant by analytical (or formal) changes of variables $\varphi$ such that $d\varphi(0)\neq 0$.
\end{proposition}

Now, with the above results, we can return to the three-dimensional world and obtain conditions for the origin of system \eqref{eq1} to be monodromic on a center manifold. We consider the following representation of \eqref{eq1}:

\begin{equation}\label{eq2}
	\begin{array}{lr}
		\dot{x}=&y+\sum\limits_{j+k+l\geqslant 2}a_{jkl}x^jy^kz^l,\\
		\dot{y}=&\sum\limits_{j+k+l\geqslant 2}b_{jkl}x^jy^kz^l,\\
		\dot{z}=&-\lambda z+\sum\limits_{j+k+l\geqslant 
			2}c_{jkl}x^jy^kz^l.
	\end{array}
\end{equation}

\begin{proposition}\label{PropoMono2}
	The origin is a nilpotent monodromic singular point with Andreev number $2$ on a center manifold of system \eqref{eq2} if and only if $b_{200}=0$ and
	\begin{equation*}%\label{ineqMono2}
		\dfrac{b_{{101}}c_{{200}}}{\lambda}<-\frac{(2a_{200}-b_{110})^2}{8}-b_{300}.
	\end{equation*}
	Moreover, if ${2a_{200}+b_{110}}\neq 0$, the restricted system satisfies the monodromy condition $\beta=n-1$ in Theorem \ref{TeoAndreevRestrito}. 
\end{proposition}

%\noindent\textbf{Proof: }The origin is monodromic only if $b_{200}=0$. Thus, by Lemma \ref{LemaMonodromiaDelta}, the necessary and sufficient condition for the Andreev number to be $2$ is $\Delta<0$, which by \eqref{eqfPhi} becomes:
%$$\frac { \left( 4\,{a_{{200}}}^{2}-4\,b_{{110}}a_{{200}}+{b_{{110}}}^{2}+8\,b_{{300}} \right) \lambda+8\,b_{{101}}c_{{200}}}{
%	\lambda}<0,$$
%or equivalently
%$$\dfrac{b_{{101}}c_{{200}}}{\lambda}<-\frac{(2a_{200}-b_{110})^2}{8}-b_{300}.$$
%%Conversely, if \eqref{ineqMono2} holds, then the hypothesis of Lemma \ref{LemaMonodromiaDelta} are satisfied.
%%
%%, a sufficient monodromy condition is given by the inequality $\Delta<0$ which becomes:
%%$$\frac { \left( 4\,{a_{{200}}}^{2}-4\,b_{{110}}a_{{200}}+{b_{{110}}}^{2}+8\,b_{{300}} \right) \lambda+8\,b_{{101}}c_{{200}}}{
%%	\lambda}<0,$$
%%equivalently 
%%$$\dfrac{b_{{101}}c_{{200}}}{\lambda}<-\frac{(2a_{200}-b_{110})^2}{8}-b_{300}.$$
%%Since under the above inequality $j^3f(0)\neq 0$ and $j^2f(0)=0$, Andreev number is $2$.
%The second statement also follows from \eqref{eqfPhi}. \qed

We highlight the singular points with Andreev number $2$ since, in the plane, there are interesting tools to study the Center Problem that only hold for this case. For instance, in \cite{GarciaFiI}, a method to solve the Center Problem for nilpotent singular points is presented which consists in finding a formal inverse integrating factor for such systems, i.e. a formal series $V(x,y)$ satisfying $XV=V\mbox{div}X$. Throughout this work, we will refer to this method as the \emph{inverse integrating factor method}. More precisely, the method is based on the following result:

\begin{theorem}[Theorem 3 in \cite{GarciaFiI}]\label{TeoGarciaFiI}
	Consider analytic system \eqref{eqPlanar} with a monodromic nilpotent singular point at the origin such that the Andreev number $n$ is even. If \eqref{eqPlanar} has a formal inverse integrating factor then the origin is a center. Moreover if $n=2$ then the converse is also true. 
\end{theorem}

The case $n=2$ was also studied in \cite{LiuLiNew1,LiuLi3Ord} where the authors also found a canonical form for nilpotent monodromic singularities of third order in the plane. In sum, the planar analytic nilpotent centers having Andreev number $2$ are completely characterized by the existence of a formal inverse integrating factor. Proposition \ref{PropoMono2} characterizes the monodromic nilpotent singular points with Andreev number $2$ in $\mathbb{R}^3$ and is a first step to extend these results to three-dimensional systems.

\subsection{Auxiliary Linear Operator}\label{SecAuxLO}

%In this section, we prove auxiliary results regarding some useful linear operators that show up frequently when working with vector fields having nilpotent singular points. These particular linear operators provide some interesting results on normal forms and methods involving formal series.

Let $H^{(n)}_{(3,1)}$ be the vector space of homogeneous polynomials of degree $n$ in three variables. 
%A basis for this vector space is given by the collection of all monomials $x^jy^kz^l$ for which $j,k,l\geqslant 0$ and $j+k+l=n$.
For a polynomial $p\in H^{(n)}_{(3,1)}$, denote by $\langle p\rangle$ be the vector subspace spanned by $p$. We consider the linear operator:
$$\begingroup
\setlength\arraycolsep{0pt}
L_n\colon\begin{array}[t]{c >{{}}c<{{}} c}
	H^{(n)}_{(3,1)} & \to & H^{(n)}_{(3,1)} \\
	\noalign{\medskip} p & \mapsto & y\dfrac{\partial p}{\partial x}-\lambda z\dfrac{\partial p}{\partial z}+\lambda p,
\end{array}
\endgroup$$
where $\lambda\neq 0$.

\begin{lemma}\label{LemaKerLnJacobi}
	The kernel of linear operator $L_n$ is given by $\ker L_n=\langle y^{n-1}z\rangle$.
\end{lemma}
%\noindent\textbf{Proof:} Trivially $L_n(y^{n-1}z)=0$. Let $p=\sum_{j+k+l=n}p_{j,k,l}x^jy^kz^l$. If $p\in\ker L_n$, then $L_n(p)=0$ which equates to:
%$$\sum_{j+k+l=n}((j+1)p_{j+1,k-1,l}-\lambda (l-1)p_{j,k,l})x^jy^kz^l=0.$$
%Consider $l\neq 1$. For $k=0$, the coefficient of $x^jz^l$ in the left-hand side of this equation is ${-\lambda (l-1)p_{j,0,l}}$ which implies that $p_{j,0,l}=0$. For $k=1$, the equation yields $p_{j,1,l}=0$. Proceeding step by step until index $k=n-1$, we conclude that $p_{j,k,l}=0$. Now, for $l=1$, the coefficient of monomial $x^jy^kz$ in the left-hand side expression above is $(j+1)p_{j+1,k-1,1}$, which implies that $p_{j,k,1}=0$ for $j>0$. Thus $p=p_{0,n-1,1}y^{n-1}z\in\langle y^{n-1}z\rangle$.\qed

\begin{lemma}\label{LemaTp+qJacobi}
	For every $q\in H^{(n)}_{(3,1)}$, there is a choice of $p\in H^{(n)}_{(3,1)}$ such that $L_n(p)+q\in\langle x^{n-1}z\rangle$.
\end{lemma}
%\noindent\textbf{Proof: }Again, we prove that $H^{(n)}_{(3,1)}=\langle x^{n-1}z\rangle\oplus \text{Im} L_n$. Lemma \ref{LemaKerLnJacobi} tells us that the codimension of Im$L_n$ is 1. As in the proof of Lemma \ref{LemaTp+q}, it is enough to prove that the intersection of $\langle x^{n-1}z\rangle$ and Im$L_n$ is $\{0\}$. We have that $L_n(p)\in\langle x^{n-1}z\rangle$ if and only if $y\dfrac{\partial p}{\partial x}-\lambda z\dfrac{\partial p}{\partial z}+\lambda p=\alpha x^{n-1}z$. The coefficient of  $x^{n-1}z$ in the left-hand side of this equation is zero, therefore $\alpha=0$. Hence $L_n(p)\in\langle x^{n-1}z\rangle$ if and only if $p\in\langle y^{n-1}z\rangle=\ker L_n$. The result holds. \qed

%\begin{remark}\label{ObsOperator}
%	Note that if we consider similar operators
%	$$\tilde{T}_n:H^{(n)}_{(3,1)}\to H^{(n)}_{(3,1)},\; p\mapsto x\dfrac{\partial p}{\partial y}-\lambda z\dfrac{\partial p}{\partial z}, \mbox{and }$$
%	$$\tilde{L}_n:H^{(n)}_{(3,1)}\to H^{(n)}_{(3,1)},\; p\mapsto y\dfrac{\partial p}{\partial x}-\lambda z\dfrac{\partial p}{\partial z}+\lambda p,$$
%	just by exchanging $x$ and $y$ in the expressions above yields analogous results. More precisely, we have $\ker\tilde{T}_n=\langle x^n\rangle$, $\ker \tilde{L}_n=\langle x^{n-1}z\rangle$  and $H^{(n)}_{(3,1)}=\langle y^n\rangle\oplus\mbox{\rm Im}\tilde{T}_n=\langle y^{n-1}z\rangle\oplus\mbox{\rm Im}\tilde{L}_n$.
%\end{remark}

\subsection{Zhitomirskii Normal Form}
Using some fundamental results in the Normal Form Theory \cite{Zhitomirskii}, we were able to obtain the following formal normal form for system \eqref{eq1}. 

\begin{theorem}[Nilpotent Normal Form in $\mathbb{R}^3$]\label{TeoFNNilpotenteZhi}
	For system \eqref{eq1} having a nilpotent singular point at the origin, there exist a formal change of variables that transforms it into the formal normal form
	\begin{equation}\label{eqN1}
		\begin{array}{lr}
			\dot{x}=y+xP_1(x),\\
			\dot{y}=Q_2(x)+yP_1(x),\\
			\dot{z}=-\lambda z+zR_1(x).
		\end{array}
	\end{equation}
	for which $P_1(0)=j^1Q_2(0)=R_1(0)=0$.
\end{theorem}

We now have the necessary tools to develop the theory in the next sections. We encourage the reader to see \cite{QueirozPessoaNillR3} for more details on the preliminary results presented up to this point. 

\section{Inverse Jacobi Multipliers and Nilpotent Centers}

\subsection{Inverse Jacobi Multipliers}\label{secJacobi}
In \cite{BuicaGarciaJacobi}, the authors explored the relationship between the center problem for Hopf points and the existence of an inverse Jacobi multiplier. Inspired by the ideas from that paper, in this section we study the properties of inverse Jacobi multipliers applied to system \eqref{eq1}, i.e. for the nilpotent case.

\begin{proposition}\label{PropoFormadoJacobiNil}
	Let $X$ be the vector field associated to system \eqref{eq1} having a nilpotent singular point in $\mathbb{R}^3$. Any non-flat $C^\infty$ inverse Jacobi multiplier $V(x,y,z)$ for vector field $X$ has the form $j^{m+1}V(0)=y^{m}z$ for $m\geqslant 0$ up to multiplication by constant.
\end{proposition}
\noindent\textbf{Proof: }Assuming that there exists an inverse Jacobi multiplier $V(x,y,z)$ for system \eqref{eq1}, it must satisfy $XV=V\mbox{div}X$. Writting $V(x,y,z)=\sum_{k=0}^{\infty}V_k$, where $V_k$ are homogeneous polynomials of degree $k$ in $x,y,z$, the following equation holds
$$\left\langle X,\sum_{k=0}^{\infty}\nabla V_k\right\rangle-\mbox{\rm div}X\left(\sum_{k=0}^{\infty}V_k\right)=0.$$

Comparing the lower order terms in both sides of the above equation yields
$$y\dfrac{\partial V_m}{\partial x}-\lambda z\dfrac{\partial V_m}{\partial z}+\lambda V_m=0,$$

where $m\geqslant 0$ is the smallest index for which $V_m\neq 0$. However the above equation is equivalent to the condition that $V_m$ is in the kernel of linear operator $L_m$ from Section \ref{SecAuxLO}. By Lemma \ref{LemaKerLnJacobi}, we must have $V_m(x,y,z)=\alpha y^{m-1}z$.\qed

The next result is proven in \cite{BuicaGarciaJacobi} and relates inverse Jacobi multipliers and invariant manifolds for generic three-dimensional systems.

\begin{theorem}[Theorem 6 in \cite{BuicaGarciaJacobi}]\label{TeoJacobiInverso}
	Let $X$ be a smooth vector field defined in an open set $U\subset\mathbb{R}^3$. Assume that $V$ is a $C^\infty$ inverse Jacobi multiplier of the form
	$$V(x,y,z)=\left(z-h(x,y)\right)W(x,y,z).$$
	Then $M=\{z=h(x,y)\}$ is an invariant manifold of $X$ and	$$v(x,y)=W(x,y,h(x,y))$$
	is an inverse integrating factor for system $X\vert_{M}$.
\end{theorem}

The following theorem is an adaptation of the result presented in \cite[Theorem 8]{BuicaGarciaJacobi} to system \eqref{eq1}.

\begin{theorem}
	Let $V(x,y,z)$ be an inverse Jacobi multiplier for vector field $X$ associated to system \eqref{eq1} having a nilpotent singular point and a local $C^\infty$ center manifold $W^c=\{z=h(x,y)\}$. Consider $V\vert_{W^c}=V(x,y,h(x,y))$. The following holds:
	\begin{itemize}
		\item[i)]$V\vert_{W^c}$ is flat at the origin;
		\item[ii)] If $W^c\subset V^{-1}(0)$, then there is a $C^\infty$-function $W(x,y,z)$ such that $W(x,y,h(x,y))\not\equiv 0$ and $V(x,y,z)=(z-h(x,y))W(x,y,z)$. Moreover $W(x,y,h(x,y))$ is an inverse integrating factor for system $X\vert_{W^c}$.
	\end{itemize}
\end{theorem}
\noindent\textbf{Proof: } First, we prove item (i). Since $W^c=\{z=h(x,y)\}$ is a $C^\infty$-parametrization for a center manifold, the following equation is satisfied:
$$\dfrac{\partial h}{\partial x}(y+P)+\dfrac{\partial h}{\partial y}Q=-\lambda h+R.$$
Let $u(x,y)=V(x,y,h(x,y))$. The above equation together with $V$ being an inverse Jacobi multiplier of $X$ yields:
\begin{eqnarray}
	\left.\dfrac{\partial u}{\partial x}(y+P)+\dfrac{\partial u}{\partial y}Q\right\vert_{z=h(x,y)}=\nonumber\\
	=\left.\dfrac{\partial V}{\partial x}(y+P)+\dfrac{\partial V}{\partial y}Q+\dfrac{\partial V}{\partial z}\left(\dfrac{\partial h}{\partial x}(y+P)+\dfrac{\partial h}{\partial y}Q\right)\right\vert_{z=h(x,y)}\nonumber\\
	=\left.\dfrac{\partial V}{\partial x}(y+P)+\dfrac{\partial V}{\partial y}Q+\dfrac{\partial V}{\partial z}(-\lambda h+R)\right\vert_{z=h(x,y)}\nonumber\\
	=\left.u(-\lambda+\mbox{div}(P,Q,R))\right\vert_{z=h(x,y)},\nonumber
\end{eqnarray}
on $W^c$. Hence, $u(0,0)=0$ since $\lambda\neq 0$. Suppose $j^\infty u(0)\neq 0$. There exists $u_m(x,y)$ a homogeneous polynomial of degree $m$ such that $j^mu(0)=u_m(x,y)$. By substituting in the above equation and comparing the lowest degree terms in both sides, we obtain $y\dfrac{\partial u_m}{\partial x}=-\lambda u_m(x,y)$. Since $u_m$ is homogeneous, it satisfies $x\dfrac{\partial u_m}{\partial x}+y\dfrac{\partial u_m}{\partial y}=m u_m$. Thus
$$\dfrac{\partial u_m}{\partial x}(\lambda x+ my)+\lambda y\dfrac{\partial u_m}{\partial y}=0.$$
Therefore $u_m(x,y)$ is a first integral for linear system $\dot{x}=\lambda x+my,\;\dot{y}=\lambda y$, which is not possible. Therefore $j^\infty u(0)=0$.

Now we prove (ii). Suppose $W^c=\{z=0\}$. By hypothesis $V(x,y,0)=0$, and by Hadamard's lemma (see \cite{DifFunctions}), $V(x,y,z)=z^kW(x,y,z)$ where $W(x,y,0)\neq 0$, $k>0$. By Proposition \ref{PropoFormadoJacobiNil}, $j^{m+1}V(0)=zy^{m}$ and hence $k=1$. The result follows.

Now if $W^c=\{z=h(x,y)\}$ with $h(x,y)\neq 0$, it is enough to consider the coordinate change $x\to x,\;y\to y,\;z\to Z=z-h(x,y)$. The new system has $Z=0$ as a center manifold and admits $\tilde{V}(x,y,Z)$ as an inverse Jacobi multiplier such that $\tilde{V}(x,y,z-h(x,y))=V(x,y,z)$.

The fact that $V(x,y,h(x,y))$ is an inverse integrating factor for $X\vert_{W^c}$ follows directly frow Theorem \ref{TeoJacobiInverso}.\qed

The search of an inverse Jacobi multiplier has its importance for the Nilpotent center problem which will be explained in Section \ref{secJacobiMeth}. The following theorem provides a direction for this search.

\begin{theorem}\label{TeoXV-divV=xnJacobi}
	Consider a vector field $X$ associated to system \eqref{eq1} having a nilpotent singular point. Then there exists a formal series $V(x,y,z)$ such that $XV-V\mbox{\rm div}X=\sum_{n\geqslant 1}\Lambda_nx^{n-1}z$.
\end{theorem}
\noindent\textbf{Proof: }Consider $V(x,y,z)=\sum_{n=1}^{\infty}V_n(x,y,z)$. Then
\begin{small}
	\begin{eqnarray}
		XV-V\mbox{\rm div}X &=&\left(y+P_2+P_3+\dots\right)\left(\dfrac{\partial V_1}{\partial x}+\dfrac{\partial V_2}{\partial x}+\dots\right)\nonumber\\
		&&+\left(-\lambda z+R_2+R_3+\dots\right)\left(\dfrac{\partial V_1}{\partial z}+\dfrac{\partial V_2}{\partial z}+\dots\right)+\nonumber\\
		&&+\left(Q_2+Q_3+\dots\right)\left(\dfrac{\partial V_1}{\partial y}+\dfrac{\partial V_2}{\partial y}+\dots\right)\nonumber\\
		&&-(V_1+V_2+\dots)\left(-\lambda+\dfrac{\partial P}{\partial x}+\dfrac{\partial Q}{\partial y}+\dfrac{\partial R}{\partial z}\right).\nonumber
	\end{eqnarray}
\end{small}
Rewriting the above expression by organizing the homogeneous terms, we have:
\begin{small}
	\begin{eqnarray}
		XV-V\mbox{div}X =\sum_{n\geqslant 1}\left(y\dfrac{\partial V_n}{\partial x}-\lambda z\dfrac{\partial V_n}{\partial z}+\lambda V_n+F_n\right)=\sum_{n\geqslant 1}L_n(V_n)+F_n,\nonumber
	\end{eqnarray}
\end{small}
where $F_n\in H^{(n)}_{(3,1)}$ are obtained by the homogeneous parts of $P,Q,R$ and $V$ of degree less than $n$. By Lemma \ref{LemaTp+qJacobi}, we can choose $V_n$ such that $L_n(V_n)+F_n=\Lambda_nx^{n-1}z$, $\Lambda_n\in\mathbb{R}$. Making the suitable choices, $XV-V\mbox{\rm div}X=\sum_{n\geqslant 1}\Lambda_nx^{n-1}z$.\qed

The next result can be found in \cite{GiacominiJacobi,GarciaFiI} and present how the existence of an inverse Jacobi multiplier is invariant by changes of variables.

\begin{theorem}\label{TeoMJIinvariant}
	Let $X$ be a vector field and $V$ is an inverse Jacobi multiplier of $X$. If there is a change of variables $\phi$ that transforms $X$ into a vector field $Y$ then $W=\det d\phi\cdot (V\circ\phi^{-1})$ is an inverse Jacobi multiplier of $Y$.	
\end{theorem}

This result can be proven by a direct, and also very long, algebraic verification. Since all the manipulations are algebraic, the above result is also true if we consider formal $\phi$ and $V$.

\subsection{Nilpotent Center Problem}\label{secJacobiMeth}

Consider planar system \eqref{eqPlanar} having a nilpotent monodromic singular point at the origin. The nilpotent center problem can be studied by investigating the existence of a formal inverse integrating factor, i.e. the \emph{inverse integrating factor method} which is based in Theorem \ref{TeoGarciaFiI}. However the existence of a formal inverse integrating factor by itself is not sufficient to prove that the singular point is a nilpotent center. The best result in this direction, also presented in \cite{GarciaFiI}, is the following.

\begin{theorem}[Theorem 4 in \cite{GarciaFiI}]%\label{TeoGarciaFiI2}
	Consider system \eqref{eqPlanar} with a monodromic nilpotent singular point at the origin having Andreev number $n$ and satisfying the condition $\beta>n-1$ given by the statement (i) in Theorem \ref{TeoAndreev}. If there is a formal inverse integrating factor $V(x,y)=\sum_{j\geqslant 2n}V_j(x,y)$ where $V_j$ are $(1,n)$-quasihomogeneous polynomials of weighted degree $j$ (see Definition \ref{defQuasihomogeneous}) and $V_{2n}\not\equiv 0$, then the origin is a center. The converse is not true.
\end{theorem}

\begin{definition}\label{defQuasihomogeneous}
	A polynomial $p\in\mathbb{R}[x,y]$ is a \emph{$(t_1,t_2)$-quasi-\-homogeneous polynomial of weighted degree $k$} if $p(\lambda^{t_1} x, \lambda^{t_2} y)=\lambda^kp(x,y)$. A general expression for such polynomials is $p(x,y)=\sum_{t_1i+t_2j=k}a_{ij}x^iy^j$ where $a_{ij}\in\mathbb{R}$. The vector space of all $(t_1,t_2)$-quasi-\-homogeneous polynomial of weighted degree $k$ is denoted by $\mathcal{P}^{(t_1,t_2)}_{k}$.
\end{definition}

We now connect the results presented in Section \ref{secJacobi} to the study of the center problem for the singular point at the origin of system \eqref{eq1} in $\mathbb{R}^3$ without restricting it to a center manifold, extending the inverse integrating factor method to the three-dimensional case. 

\begin{definition}
	Let $X$ be a formal vector field, i.e. $X=X_1\partial_x+X_2\partial_y+X_3\partial_z$ where $X_i\in \mathbb{R}[[x,y,z]]$ for $i=1,2,3$ and $M\in \mathbb{R}[[x,y,z]]$. We say that $M=0$ is a \emph{formal center manifold} for $X$ when $M$ formally satisfies the equation $XM=K\cdot M$ for some $K\in \mathbb{R}[[x,y,z]]$ and $\nabla M (0)\wedge (0,0,1)^T=0$.
\end{definition}

Note that when $X,M,K$ are analytical, $M=0$ is an analytical center manifold for $X$. Let $(u,v,w)=\phi(x,y,z)$ be a formal near-identity change of variables, i.e. $d\phi(0)=Id$. We denote the vector field in the new variables by $Y$, and it satisfies $d\phi\cdot X=Y\circ\phi$. Thus, if $M=0$ is a formal center manifold for $X$, then, for $L=M\circ\phi^{-1}$, $L=0$ is a formal center manifold for $Y$. In fact, we have $XM=K\cdot M$ and, by definition:
\begin{eqnarray}
YL&=&\langle \nabla L,Y\rangle=\langle\nabla M d\phi^{-1}, d\phi X\circ\phi^{-1}\rangle=\langle \nabla M,X\rangle\circ \phi^{-1}=\nonumber\\
&=&(K\cdot M)\circ\phi^{-1}=(K\circ\phi^{-1})\cdot L.\nonumber
\end{eqnarray}

Moreover, $\nabla L(0)=\nabla M(0) d\phi^{-1}(0)=\nabla M(0)$. For $\bar{\phi}=\phi\vert_{M=0}$, we have
$$d\bar{\phi}\cdot X\vert_{M=0}=(d\phi\cdot X)\vert_{M=0}=(Y\circ\phi)\vert_{M=0}=Y\vert_{L=0}\circ\bar{\phi}.$$
Thus, $\bar{\phi}$ is a near-identity formal change of variables that transforms $X\vert_{M=0}$ into $Y\vert_{L=0}$. Therefore, by Theorem \ref{TeoMJIinvariant}, the existence of an inverse integrating factor for the restriction $X\vert_{M=0}$ would imply the existence of an inverse integrating factor for any system $Y=d\phi X\circ\phi^{-1}$ restricted to its formal center manifold $L=0$.

%Note that since $\nabla L_0\wedge(0,0,1)^T=0$, $L=0$ admits a formal parametrization $w-g(u,v)=0$ \cite[Proposition 3.1]{SokalImplicit} where $g\in\mathbb{R}[[u,v]]$.
%
%Now consider vector field $X$ associated to system \eqref{eq1} having a nilpotent singular point at the origin. Solving recursevily equation \eqref{eqAlgoWc} it is possible to obtain an expression for a formal center manifold $z-h(x,y)=0$ which we will denote by $W^c$. The existence of an inverse integrating factor for the restriction $X\vert_{W^c}$ would imply the existence of an inverse integrating factor for any system $Y=d\phi X\circ\phi^{-1}$ restricted to its formal center manifold $L=0$.

Now consider vector field $X$ associated to system \eqref{eq1} having a nilpotent singular point at the origin. By Theorem \ref{TeoFNNilpotenteZhi}, there is a formal near-identity change of variables $\phi$ such that $Y$ is given by:
\begin{equation*}
	\begin{array}{lr}
		\dot{x}=y+xP_1(x)=Y_1,\\
		\dot{y}=Q_2(x)+yP_1(x)=Y_2,\\
		\dot{z}=-\lambda z+zR_1(x)=Y_3.
	\end{array}
\end{equation*}
Since $Y$ is decoupled, its restriction to $L=0$ is given by the first two components. Consequentially, if $\tilde{v}(x,y)$ is an inverse integrating factor for $Y\vert_{L=0}$, then $z\tilde{v}(x,y)$ is an inverse Jacobi multiplier for $Y$. In fact:
\begin{eqnarray}
Y(z\tilde{v})-z\tilde{v}\mbox{\rm div}Y&=zY_1\dfrac{\partial\tilde{v}}{\partial x}+zY_2\dfrac{\partial\tilde{v}}{\partial y}+Y_3\tilde{v}-z\tilde{v}\left(\mbox{\rm div}(Y_1,Y_2)+\dfrac{\partial Y_3}{\partial z}\right)=\nonumber\\
&=\tilde{v}\left(Y_3-z\dfrac{\partial Y_3}{\partial z}\right)=0.\nonumber
\end{eqnarray}
Therefore, by Theorem \ref{TeoMJIinvariant} the original vector field $X$ would also admit an inverse Jacobi multiplier.

Using this argument, we can prove Theorem \ref{TeoJacobin=2}.
\medskip
%
%\begin{theorem}\label{TeoJacobin=2}
%	Consider system \eqref{eq1} having a monodromic singular point with Andreev number $2$ at the origin. If the origin is an analytic nilpotent center on a center manifold, then there exists a formal inverse Jacobi multiplier $V(x,y,z)$ for system \eqref{eq1} such that $j^{m+1}V(0)=y^mz$ for $m\geqslant 0$.
%\end{theorem}

\noindent\textbf{Proof of Theorem \ref{TeoJacobin=2}: }Suppose the origin is an analytic nilpotent center on a center manifold of vector field $X$ associated to system \eqref{eq1}. Denote by $\bar{X}$ the analytic restricted vector field to the center manifold. Since the Andreev number for the origin of the restricted system is $2$, by Theorem \ref{TeoGarciaFiI} there exists a formal series $v(x,y)$ satisfying $\bar{X}v-v\mbox{\rm div}\bar{X}=0$, i.e. $v(x,y)$ is a formal inverse integrating factor for $\bar{X}$. Therefore, there exists an inverse integrating factor $\tilde{v}(x,y)$ for the vector field associated to the first two components of the normal form \eqref{eqN1}. Thus, normal form \eqref{eqN1} admits inverse Jacobi multiplier $\tilde{V}(x,y,z)=z\tilde{v}(x,y)$ which implies that the original system \eqref{eq1} also has a formal inverse Jacobi multiplier $V(x,y,z)$. By Proposition \ref{PropoFormadoJacobiNil}, we must have $j^{m+1}V(0)=zy^{m}$, for some $m\geqslant 0$. The result follows.\qed

\begin{remark}
	Since the formal change of variables that transforms system \eqref{eq1} into normal form \eqref{eqN1} is a near-identity one, the first non-zero jet of any formal inverse Jacobi multiplier is preserved.
\end{remark}

The idea is to use the above results to devise an algorithm to solve the nilpotent center problem for systems \eqref{eq1} with Andreev number $2$. We know, by Proposition \ref{PropoFormadoJacobiNil}, that any inverse Jacobi multiplier of system \eqref{eq1} must have $j^{m+1}V(0)=zy^m$. However, without reducing the possible values for $m$, such algorithm would not be practical for computational purposes. The next results, proven in \cite{GarciaFiI}, determine the possible values for $m$.

\begin{proposition}\label{PropoGarciaFormaFiIbeta=n-1}
	Consider system \eqref{eqPlanar} having monodromic singular point at the origin with Andreev number $n$, satisfying monodromy condition $\beta=n-1$. If there is a formal inverse integrating factor $V(x,y)$ for \eqref{eqPlanar}, then it has the form $j^2V(0)=y^2$.
\end{proposition}

\begin{proposition}\label{PropoGarciaFormaFiIbeta>n-1}
	Consider system \eqref{eqPlanar} having monodromic singular point at the origin with Andreev number $n$, satisfying monodromy condition $\beta>n-1$ in Theorem \ref{TeoAndreev}. If the system admits formal inverse integrating factor $V(x,y)$, then either $V(0,0)\neq 0$ or $j^2V(0)=y^2$.
\end{proposition}

Thus, for the three-dimensional case, any inverse Jacobi multiplier have either $j^1V(0)=z$ or $j^3V(0)=zy^2$.

\begin{remark}
	Propositions \ref{PropoGarciaFormaFiIbeta=n-1} and \ref{PropoGarciaFormaFiIbeta>n-1} are direct consequences of Propositions 7 and 8 in \cite{GarciaFiI}, respectively, since the monodromy conditions are invariant under changes of variables (see Proposition \ref{Propobeta=n-1invariante}).%\cite{QueirozPessoaNillR3}
\end{remark}

Now, we can use the above results to apply the following procedure to find necessary conditions for the origin of system \eqref{eq1} with Andreev number $2$ to be an analytic nilpotent center. More precisely, we use Proposition \ref{PropoMono2} to determine which of the monodromy conditions in Theorem \ref{TeoAndreevRestrito} are satisfied. Then, we consider a formal series $V(x,y,z)=\sum_{n\geqslant 1}V_n(x,y,z)$, where $V_n\in H^{(n)}_{3,1}$. By Theorem \ref{TeoXV-divV=xnJacobi}, it is always possible to find $V(x,y,z)$ such that
$$XV-V\mbox{\rm div}X=\sum_{n\geqslant 1}\Lambda_n x^{n-1}z.$$
The quantities $\Lambda_n$ are polynomials in the parameters of system \eqref{eq1}. Any non-zero $\Lambda_n$ is an obstruction for $V(x,y,z)$ to be an inverse Jacobi multiplier. If there are obstructions for both possibilities $j^1V(0)=z$ and $j^3V(0)=zy^2$ then the origin cannot be an analytic nilpotent center on the center manifold. We apply this algorithm in Section \ref{SecExamples}.

\section{Analytic Nilpotent Centers as limits of Non-degenerated Centers in $\mathbb{R}^3$}\label{SecLlimite}

A planar vector field $X$ has a non-hyperbolic non-degenerate singular point $p$ if the eigenvalues of $DX(p)$ are purely imaginary and conjugated, say $\pm\beta i$ for $\beta\neq 0$. One canonical form for such systems is
\begin{equation}\label{eqPlanarNonDeg}
	\begin{array}{lr}
		\dot{x}=y+P(x,y),\\
		\dot{y}=-x+Q(x,y),
	\end{array}
\end{equation}
where $j^1P(0)=j^1Q(0)=0$. This singular point is monodromic. In regards to the Center Problem, the following result, due to Poincaré and Lyapunov, holds true:

\begin{theorem}[Poincaré-Lyapunov Center Theorem]\label{TeoPoincareLyapunov}
	The origin of system \eqref{eqPlanarNonDeg} is a center if and only if there is an analytical first integral $H(x,y)$ such that $j^2H(0)=x^2+y^2$.
\end{theorem}

A proof for this result can be found in \cite{Poincare,Queiroz}. The Poincaré-Lypunov Center Theorem provides an algorith to solve the center problem which consists in the search for a first integral. Such algoritm, refered to as \emph{Poincaré-Lyapunov algorithm}, is well-known and extensively studied in the literature (see, for instance \cite{AndronovB,Llibre,Queiroz,Romanovski}). 

For the nilpotent case, the authors of \cite{LlibreLimite} (and its previous versions) were able to use the Poincaré-Lyapunov algorithm to investigate the nilpotent center problem for planar systems \eqref{eqPlanar}. The key of their approach is the following theorem.

\begin{theorem}\label{TeoLlibreLimite}
	Suppose the origin of the following system
	\begin{equation*}
		\begin{array}{lr}
			\dot{x}=y+P(x,y),\\
			\dot{y}=Q(x,y),
		\end{array}
	\end{equation*}
	is a nilpotent center. Then there are two functions $F_1(x,y)$, $F_2(x,y)$ analytical at the origin with $j^1F_1(0)=j^1F_2(0)=0$ such that the 1-parameter family
	\begin{equation*}
		\begin{array}{lr}
			\dot{x}=y+P(x,y)+\varepsilon F_1(x,y),\\
			\dot{y}=-\varepsilon x+ Q(x,y)+\varepsilon F_2(x,y).
		\end{array}
	\end{equation*}
	has a non-degenerate center at the origin for any $\varepsilon>0$. Also there is an analytic function $f(x,y)$ at the origin with $j^1f(0)=0$ such that $(x-F_1(x,y))\dfrac{\partial f}{\partial y}=F_2(x,y)(1+\dfrac{\partial f}{\partial x})$.
\end{theorem}

The above result tells us that for the analytical planar case, any nilpotent center is limit of non-degenerate centers. We now proceed to study a possible extension of this result to $\mathbb{R}^3$. 

A three-dimensional system has a Hopf singular point if it can be put, via changes of variables and time rescalings, in the form:
\begin{equation}\label{eqNonDeg}
	\begin{array}{lr}
		\dot{x}=y+P(x,y,z),\\
		\dot{y}=-x+Q(x,y,z),\\
		\dot{z}=-\lambda z+R(x,y,z),
	\end{array}
\end{equation}
for which $j^1P(0)=j^1Q(0)=j^1R(0)=0$ and $\lambda\neq 0$. That is, the eigenvalues of its correspondent Jacobian matrix are two pure imaginary numbers and one non-zero real number. The Center Manifold Theorem also holds for system \eqref{eqNonDeg} and thus it is also possible to study the Center Problem for Hopf singular points. For this matter, Theorem \ref{TeoPoincareLyapunov} has its three-dimensional counterpart (a proof is provided in \cite{Bibikov,Queiroz}): 

\begin{theorem}[Lyapunov Center Theorem]%\label{TeoLyapunovR3}
	The origin of system \eqref{eqNonDeg} is a center on a center manifold if and only if there is an analytical first integral $H(x,y,z)$ such that $j^2H(0)=x^2+y^2$. Moreover, if the origin is a center on a center manifold, the center manifold is unique and analytical.
\end{theorem}

There are several obstacles for the extension of Theorem \ref{TeoLlibreLimite} to three-dimensional system. Among them, the lack of analiticity of the center manifold is one of the most noticeable. Therefore we make some assumptions to obtain a first result. A more detailed study on the center manifolds can be found in \cite{Kelley, Sijbrand}.

%\begin{theorem}\label{TeoLlibreLimiteR3}
%	Suppose the origin of the following system
%	\begin{equation}\label{eq1repete}
%		\begin{array}{lr}
%			\dot{x}=y+P(x,y,z),\\
%			\dot{y}=Q(x,y,z),\\
%			\dot{z}=-\lambda z+R(x,y,z),
%		\end{array}
%	\end{equation}
%	is a nilpotent center on an analytical center manifold. Then there are two functions $F_1(x,y)$, $F_2(x,y)$ analytical at the origin with $j^1F_1(0)=j^1F_2(0)=0$ such that the 1-parameter family
%	\begin{equation}\label{eq1Llimit}
%		\begin{array}{lr}
%			\dot{x}=y+P(x,y,z)+\varepsilon F_1(x,y),\\
%			\dot{y}=-\varepsilon x+Q(x,y,z)+\varepsilon F_2(x,y),\\
%			\dot{z}=-\lambda z + R(x,y,z).
%		\end{array}
%	\end{equation}
%	has a non-degenerate center at the origin for any $\varepsilon>0$. Also there is an analytic function $f(x,y)$ at the origin with $j^1f(0)=0$ such that $(x-F_1(x,y))\dfrac{\partial f}{\partial y}=F_2(x,y)(1+\dfrac{\partial f}{\partial x})$.
%\end{theorem}
\medskip
\noindent\textbf{Proof of Theorem \ref{TeoLlibreLimiteR3}: }By hypothesis, \eqref{eq1} admits a center manifold $z=h(x,y)$ on which the origin is an analytic nilpotent center. Thus, the restricted system is analytical and has a nilpotent center at the origin. Consider the following change of variables 
$$x=x,\;y=y,\;Z=z-h(x,y).$$
%Note that $\dot{Z}=\dot{z}-\dfrac{\partial h}{\partial x}\dot{x}-\dfrac{\partial h}{\partial y}\dot{y}=-\lambda Z+h(x,y)+R(x,y,Z+h(x,y))-\dfrac{\partial h}{\partial x}\dot{x}-\dfrac{\partial h}{\partial y}\dot{y}=-\lambda Z+\tilde{R}(x,y,Z),$ onde $\tilde{R}(x,y,0)=0$. Assim, nas novas coordenadas, o sistema se escreve como:
Under these new coordinates, system \eqref{eq1} becomes
\begin{equation*}
	\begin{array}{lr}
		\dot{x}=y+P(x,y,Z+h(x,y)),\\
		\dot{y}=Q(x,y,Z+h(x,y)),\\
		\dot{Z}=-\lambda Z + Z\tilde{R}(x,y,Z).
	\end{array}
\end{equation*}
We have that $Z=0$ is an analytical center manifold. Since the restricted system is analytical, by Theorem \ref{TeoLlibreLimite} there exist two functions $F_1(x,y),F_2(x,y)$ analytical at the origin with $j^1F_1(0)=j^1F_2(0)=0$ such that 
\begin{equation*}
	\begin{array}{lr}
		\dot{x}=y+P(x,y,h(x,y))+\varepsilon F_1(x,y),\\
		\dot{y}=-\varepsilon x+ Q(x,y,h(x,y))+\varepsilon F_2(x,y)
	\end{array}
\end{equation*}
has a non-degenerate center at the origin for every $\varepsilon>0$. Hence, the three-dimensional system
\begin{equation*}
	\begin{array}{lr}
		\dot{x}=y+P(x,y,Z+h(x,y))+\varepsilon F_1(x,y),\\
		\dot{y}=-\varepsilon x+Q(x,y,Z+h(x,y))+\varepsilon F_2(x,y),\\
		\dot{z}=-\lambda Z + Z\tilde{R}(x,y,Z).
	\end{array}
\end{equation*}
has a non-degenerate center on $Z=0$. Returning to the original coordinates, we conclude that system
\begin{equation*}
	\begin{array}{lr}
		\dot{x}=y+P(x,y,z)+\varepsilon F_1(x,y),\\
		\dot{y}=-\varepsilon x+Q(x,y,z)+\varepsilon F_2(x,y),\\
		\dot{z}=-\lambda z + R(x,y,z).
	\end{array}
\end{equation*}
has a non-degenerate center at the origin on center manifold $z=h(x,y)$. \qed

\begin{remark}The equation 
	$$(x-F_1(x,y))\dfrac{\partial f}{\partial y}=F_2(x,y)(1+\dfrac{\partial f}{\partial x}),$$
	implies that $x$ must factor out the lowest degree homogeneous polynomial of the power series expansion of $F_1$. Thus, the coefficient of $y^2$ in its expansion must be zero.
\end{remark}

\begin{remark}
	The arguments presented in the previous proof are also valid even when the center manifold $z=h(x,y)$ is not analytical as long as the restricted system is. For instance, the decoupled systems given by
	\begin{equation*}
		\begin{array}{lr}
			\dot{x}=y+P(x,y),\\
			\dot{y}=Q(x,y),\\
			\dot{z}=-\lambda z + R(x,y,z).
		\end{array}
	\end{equation*}
	have this property since its restriction to any center manifold is completely determined by the first two components.
\end{remark}

Now, we use the Poincaré-Lyapunov method to detect  analytic nilpotent centers. First, we consider two formal series: 
$$F_i=\sum_{j+k\geqslant 2}g_{jk}^{i}x^jy^k,\;i=1,2,$$
with $g_{02}^{1}=0$. We then perturb system \eqref{eq1} as in \eqref{eq1Llimit} and search for obstructions for the formal series $H(x,y,z)=\varepsilon x^2+y^2+\sum_{n\geqslant 3}H_n(x,y,z)$ to be a first integral. Those obstructions will be the \emph{Lyapunov quantities}. More precisely, let $X_\varepsilon$ be the vector field associated to \eqref{eq1Llimit}. We can always find $H(x,y,z)$ such that
$$X_\varepsilon H=\sum_{l\geqslant 1}\eta_{l}(x^2+y^2)^l,$$
where $\eta_{l}$ is a rational function on the parameters of system \eqref{eq1Llimit}, i.e. the coefficients of functions $P,Q,R$, the coefficients $g_{jk}^{i}$ as well as $\varepsilon$ and $\lambda$. The numerator of the Lyapunov quantities will be a polynomial on $\varepsilon$. For the system \eqref{eq1Llimit} to have a center on the center manifold for every $\varepsilon>0$, we must have all the Lyapunov quantities to be null. We illustrate this method in the examples of the next section.

\section{Applications}\label{SecExamples}

\subsection{The Song-Wang-Feng system} The following system was first considered in \cite{SongWangFeng}.
\begin{equation}\label{eqSongFengWang}
	\begin{array}{lcr}
		\dot{x}=y-2xy+axz,\\
		\dot{y}=-2x^3+y^2+byz,\\
		\dot{z}=-z+dxy.
	\end{array}
\end{equation}
In \cite{SongWangFeng} the authors solved the center problem by considering a polynomial approximation of the center manifold. Here, we solve the center problem without restricting the considered system to a center manifold. 

By Proposition \ref{PropoMono2}, the origin of \eqref{eqSongFengWang} is monodromic with Andreev number $n=2$ and satisfies the monodromy condition $\beta>n-1$.

Using the algorithm described in Section \ref{SecLlimite}, we compute the Lyapunov quantities of the perturbation \eqref{eq1Llimit} for system \eqref{eqSongFengWang}. We have:
$$\eta_{1}=-\dfrac{4\varepsilon^2d(a-b)}{12\varepsilon+3}.$$
By Theorem \ref{TeoLlibreLimiteR3}, $d(a-b)=0$ is a necessary condition for the origin to be an analytic nilpotent center.

Under $d=0$, system \eqref{eqSongFengWang} becomes
\begin{equation*}
	\begin{array}{lcr}
		\dot{x}=y-2xy+axz,\\
		\dot{y}=-2x^3+y^2+byz,\\
		\dot{z}=-z.
	\end{array}
\end{equation*}
which has $z=0$ as an analytical center manifold. The restricted system is the Hamiltonian system given by
\begin{equation}\label{eqSongFengWangRestricted}
	\begin{array}{lcr}
		\dot{x}=y-2xy,\\
		\dot{y}=-2x^3+y^2.
	\end{array}
\end{equation}
Thus, $d=0$ implies that the origin is a nilpotent center on a center manifold. Now we search for obstructions for system \eqref{eqSongFengWang} to have a formal inverse Jacobi multiplier. We consider both possibilities of $V(x,y,z)$. First, if $j^1V(0)=z$. We have $\Lambda_1=\dots=\Lambda_4=0$ and:
$$\Lambda_5=-4d\,(2a-b),\;\Lambda_6=-a\,d.$$ 
Now for $j^3V(0)=zy^2$, we obtain $\Lambda_1=\dots=\Lambda_8=0$ and:
$$\Lambda_9=-{\frac {12\,d \left( 2\,a-b \right) }{5}},\;\Lambda_{10}=-\dfrac{2\,d \left( 9\,a-2\,b \right)}{15}.$$ Thus, by Theorem \ref{TeoJacobin=2}, the conditions $a=b=0$ or $d=0$ are necessary for the origin to be an analytic nilpotent center on a center manifold. 

Now, $a=b=0$ is a center condition, since, under this condition, the first two components of \eqref{eqSongFengWang} do not depend on $z$. Thus, the restricted system is also given by \eqref{eqSongFengWangRestricted} regardless the center manifold. Therefore $a=b=0$ or $d=0$ are sufficient conditions for system \eqref{eqSongFengWang} to have an analytic nilpotent center on a center manifold. We have proven the following result:

\begin{theorem}%\label{TeoSongWangFeng}
The origin of system \eqref{eqSongFengWang} is an analytic nilpotent center on a center manifold if and only if $a=b=0$ or $d=0$.
\end{theorem}

\subsection{Jerk system}
Consider the following jerk system
\begin{equation}\label{eqJerk3}
	\begin{array}{lcr}
		\dot{x}=y-F(x,y,z),\\
		\dot{y}=F(x,y,z)\\
		\dot{z}=-z+F(x,y,z),
	\end{array}
\end{equation}
where $F(x,y,z)=g_{300}x^3+g_{011}yz+g_{210}x^2y+g_{120}xy^2+g_{030}y^3$. Suppose that the origin has Andreev number $n=2$, which by Proposition \ref{PropoMono2} implies that $g_{300}>0$. Note that the monodromy condition $\beta>n-1$ is satisfied. Applying Theorem \ref{TeoXV-divV=xnJacobi} with the formal series $V(x,y,z)=\sum_{j+k+l\geqslant 1}v_{jkl}x^jy^kz^l$, we compute the obstructions $\Lambda_n$ for the existence of an Inverse Jacobi Multiplier. For $j^1V(0)=z$, we have $\Lambda_1=\Lambda_2=0$ and:
\begin{eqnarray}
\Lambda_3&=&3g_{{300}}-g_{{210}},\nonumber\\
\Lambda_4&=&-g_{{300}} \left( g_{{011}}-v_{{011}} \right),\nonumber\\
\Lambda_5&=&3\,g_{{300}} \left( g_{{030}}-g_{{120}}+2\,g_{{300}} \right)\;\mod\langle \Lambda_3,\Lambda_4\rangle,\nonumber\\
\Lambda_6&=&g_{{300}}g_{{011}} \left( g_{{120}}+3\,g_{{300}} \right)\;\mod\langle \Lambda_3,\Lambda_4,\Lambda_5\rangle,\nonumber\\
\Lambda_7&=&-4\,{g_{{300}}}^{2}{g_{{011}}}^{2}\;\mod\langle\Lambda_3,\Lambda_4,\Lambda_5\rangle.\nonumber
\end{eqnarray}
And for $j^3V(0)=zy^2$, we obtain $\Lambda_1=\dots=\Lambda_6=0$ and:
\begin{eqnarray}
	\Lambda_7&=&\tfrac{1}{6}\,g_{{300}} \left(g_{{210}}-3\,g_{{300}} \right), \nonumber\\
	\Lambda_8&=&\tfrac{1}{4}{g_{{300}}}^{2} \left( g_{{011}}-3\,v_{{031}} \right), \nonumber\\
	\Lambda_9&=&-{\frac {9\,{g_{{300}}}^{2} \left( g_{{030}}+2\,g_{{300}}-g_{{120}} \right) }{10}}\;\mod\langle\Lambda_7,\Lambda_8\rangle\nonumber\\
	\Lambda_{10}&=&-\tfrac{1}{3}g_{{011}}{g_{{300}}}^{2} \left( 3\,g_{{300}}+g_{{120}}
	\right)\;\mod\langle\Lambda_7,\Lambda_8,\Lambda_9\rangle,\nonumber\\
	\Lambda_{11}&=&{\frac {10\,{g_{{011}}}^{2}{g_{{300}}}^{3}}{7}}\;\mod\langle\Lambda_7,\Lambda_8,\Lambda_9\rangle\nonumber.	
\end{eqnarray}
Both the possibilities above yield the same conditions on the parameters for the existence of an inverse Jacobi multiplier: $g_{030}=g_{120}-2g_{300}$, $g_{011}=0$ and $g_{210}=3g_{300}$. Thus, by Theorem \ref{TeoJacobin=2}, the origin can only be an analytic center on a center manifold if those conditions are satisfied. These are in fact center conditions, since under those, the restricted system is given by the first two equations of \eqref{eqJerk3} and has
\begin{eqnarray}
v(x,y)=1+\left( g_{{120}}-3\,g_{{300}} \right){x}^{2}  +2 \left( g_{{120}}-3\,g_{{300}} \right)xy\nonumber\\
+{\frac {\left( g_{{120}}-3\,g_
		{{300}} \right)  \left( g_{{120}}-2\,g_{{300}} \right){y}^{2} }{g_{{300}}}},\nonumber
\end{eqnarray}
as an inverse integrating factor. By Theorem \ref{TeoGarciaFiI}, the origin is a center. We conclude:

\begin{theorem}
The origin of system \eqref{eqJerk3} is a nilpotent center on a center manifold if and only if $g_{030}=g_{120}-2g_{300}$, $g_{011}=0$ and $g_{210}=3g_{300}$.
\end{theorem}

\subsection{A system with fixed Andreev number}
If a planar system has a monodromic nilpotent singular point with Andreev number $n$, it can be put in the following form:
\begin{equation*}
	\begin{array}{lcr}
		\dot{x}=y+\mu x^n+\sum_{k\geqslant n+1}P_k(x,y),\\
		\bigskip
		\dot{y}=-nx^{2n-1}+n\mu x^{n-1}y+\sum_{k\geqslant 2n}Q_k(x,y),\\
	\end{array}
\end{equation*}
where $P_k,Q_k\in\mathcal{P}^{(1,n)}_k$ (see \cite{GarciaGineInte}, for more details). 
Thus, an interesting family of three-dimensional systems to consider are the following
\begin{equation}\label{eqAndreevfixo}
	\begin{array}{lcr}
		\dot{x}=y+\mu x^n+\sum_{k\geqslant n}P_k(x,y,z),\\
		\dot{y}=-nx^{2n-1}+n\mu x^{n-1}y+\sum_{k\geqslant 2n-1}Q_k(x,y,z),\\
		\dot{z}=-z+R(x,y,z),
	\end{array}
\end{equation}
where $P_k, Q_k$ are $(1,n,1)$-quasi-homogeneous polynomials of weighted degree $k$ and $j^1R(0)=0$. For any parametrization $z=h(x,y)$ of a center manifold, we have that the functions $F(x),f(x)$ and $\Phi(x)$ in Theorem \ref{TeoAndreevRestrito} are given by:
$$F(x)=-\mu x^n+O(x^{n+1}),$$
$$f(x)=-n(1+\mu^2)x^{2n-1}+O(x^{2n}),$$
$$\Phi(x)=2n\mu x^{n-1}+O(x^n).$$
Thus, we have that the origin is monodromic with Andreev number $n$. Moreover, the parameter $\mu$ defines whether we have monodromy conditions $\beta=n-1$ or $\beta>n-1$. 

We solve the nilpotent center problem for a particular case of system \eqref{eqAndreevfixo} with $n=2$:
\begin{equation}\label{eqAndreevfixo2}
	\begin{array}{lcr}
		\dot{x}=y+\mu x^2+a_{101}xz+a_{002}z^2,\\
		\dot{y}=-2x^{3}+2\mu xy,\\
		\dot{z}=-z+c_{200}x^2+c_{110}xy+c_{020}y^2.
	\end{array}
\end{equation}
%Using Theorem \ref{TeoXH=xn}, we compute the quantities $\omega_n$.  We have that $\omega_{5}=8\mu(1+\mu^2)$. We begin searching for integrable centers, i.e. the case $\mu=0$. We compute the quantities $\omega_{6},\dots,\omega_{13}$. The parameter values for which all these quantities vanish are $\mu=0$ and either $a_{101}=a_{002}=0$ or $c_{200}=c_{110}=c_{020}=0$. For both possibilities, the restricted system is the Hamiltonian system $\dot{x}=y$, $\dot{y}=-2x^3$. These are the only integrable nilpotent centers for system \eqref{eqAndreevfixo2}.

Using Theorem \ref{TeoXV-divV=xnJacobi}, we compute the obstructions for the existence of an Inverse Jacobi Multiplier $V(x,y,z)$. For the case $j^1V(0)=z$, we compute the obstructions $\Lambda_1,\dots,\Lambda_{10}$. Due to the size of these quantities, we exhibit only $\Lambda_1,\dots,\Lambda_5$ and omit $\Lambda_{6},\dots,\Lambda_{10}$.
\begin{eqnarray}
\Lambda_1&=&0,\;\Lambda_2=-4\mu,\nonumber\\
\Lambda_3&=&-3\,a_{{101}}c_{{200}},\nonumber\\
\Lambda_4&=&-2\,{\mu}^{2}v_{{011}}-2\,\mu\,a_{{101}}c_{{110}}+4\,\mu\,a_{{101}}c_{{200}}-4\,a_{{002}}{c_{{200}}}^{2}-2\,v_{{011}},\nonumber\\
\Lambda_5&=&-\tfrac{1}{3}\,c_{{020}}{\mu}^{2}a_{{101}}+2\,c_{{110}}{\mu}^{2}a_{{101}}-24\,c_{{200}}{\mu}^{2}a_{{101}}-2\,c_{{110}}a_{{002}}\mu\,c
_{{200}}\nonumber\\
&&-2\,c_{{200}}\mu\,a_{{101}}v_{{011}}-c_{{110}}c_{{200}}{a_{{101}}}^{2}-8\,a_{{101}}c_{{110}}+16\,a_{{101}}c_{{200}}\nonumber\\
&&+3\,{c_{{200}}}^{2}{a_{{101}}}^{2}-2\,c_{{020}}a_{{101}}.\nonumber
\end{eqnarray}

Setting $\Lambda_2=\dots=\Lambda_{10}=0$ yields the conditions on the parameters: $\mu=0$ and $a_{101}=a_{002}=0$ or $c_{200}=c_{110}=c_{020}$. 

For the case $j^3V(0)=zy^2$, we compute the obstructions $\Lambda_1,\dots,\Lambda_{15}$. Due to the size of these quantities, we exhibit only $\Lambda_1,\dots,\Lambda_9$ and omit $\Lambda_{10},\dots,\Lambda_{15}$. We have $\Lambda_1=\dots=\Lambda_6=0$ and 

\begin{eqnarray}
\Lambda_7&=&-\dfrac{c_{{200}}a_{{101}} \left( {\mu}^{2}+5 \right)}{5},\nonumber\\
\Lambda_8&=&\tfrac{\mu}{3}
\left( {\mu}^{2}+5 \right) a_{{101}}c_{{110}}+\tfrac{2\mu}{3}\left( {
	\mu}^{2}+1 \right) c_{{200}}a_{{101}}\nonumber\\
&&-\tfrac{2}{3} \left( {\mu}^{2}+3
\right) a_{{002}}{c_{{200}}}^{2}-2/5\,\mu\,{a_{{101}}}^{2}{c_{{
			200}}}^{2}-\tfrac{1}{3} \left( {\mu}^{4}+10\,{\mu}^{2}+9 \right) v_{{031}},\nonumber\\
\Lambda_9&=&-{\frac { \left( 54\,{\mu}^{4}+222\,{\mu}^{2}+168 \right) a_{{101}}c
		_{{110}}}{35}}+{\frac { \left( 48\,{\mu}^{2}+48 \right) c_{{200}}a_{{101}}}{5}}\nonumber\\
	&&-{\frac { \left( 27\,{\mu}^{4}+129\,{\mu}^{2}+42
		\right) a_{{101}}c_{{020}}}{35}}+\tfrac{\mu}{35}\left( 11\,{\mu}^{2}
-29 \right) c_{{200}}a_{{101}}v_{{031}}\nonumber\\
&&+{\frac {72\,\mu\,
		\left( {\mu}^{2}+1 \right) a_{{002}}{c_{{200}}}^{2}}{35}}+{\frac 
	{6\,\mu\, \left( 9\,{\mu}^{2}+29 \right) c_{{200}}a_{{002}}c_{{110}}}{35}}\nonumber\\
&&+\tfrac{1}{35}\left( 34\,{\mu}^{2}+35 \right) {a_{{101}}}^{2}c_{
	{110}}c_{{200}}+{\frac { \left( 57\,{\mu}^{2}-105 \right) {a_{{101}}}^{2}{c_{{200}}}^{2}}{35}}\nonumber\\
&&-\tfrac{6\mu}{5}\,a_{{101}}\,a_{{002}}{c_{{200}}}^{3}.\nonumber
\end{eqnarray}

Setting $\Lambda_7=\dots=\Lambda_{15}=0$ yields the following conditions on the parameters: $a_{101}=a_{002}=0$ or $c_{200}=c_{110}=c_{020}=0$. 

The conditions for \eqref{eqAndreevfixo2} to have any inverse Jacobi multiplier are $a_{101}=a_{002}=0$ or $c_{200}=c_{110}=c_{020}=0$. By Theorem \ref{TeoJacobin=2}, those are necessary conditions for the origin to be an analytic nilpotent center.

These conditions are also sufficient since under either one of those, the restricted system is 
\begin{equation*}
	\begin{array}{lcr}
		\dot{x}=y+\mu x^2,\\
		\dot{y}=-2x^{3}+2\mu xy.
	\end{array}
\end{equation*}
which is invariant by the transformations $x\mapsto -x$, $y\mapsto y$, $t\mapsto -t$.

%\begin{figure}[h]
%	\centering
%	\subfigure[$\mu<0.$]{
%		\includegraphics[scale=0.4]{Andreevfixo2muneg}}
%	\subfigure[$\mu=0$ (Hamiltonian).]{
%		\includegraphics[scale=0.4]{Andreevfixo2mu0}}
%	\subfigure[$\mu>0.$]{
%		\includegraphics[scale=0.4]{Andreevfixo2mupos}}
%	\caption{Phase portraits of the restriction systems \eqref{eqAndreevfixo2} under $a_{101}=a_{002}=0$ or $c_{200}=c_{110}=c_{020}=0$.}
%\end{figure}

Therefore we conclude:
\begin{theorem}
	The origin of system \eqref{eqAndreevfixo2} is an analytic nilpotent center on a center manifold if and only if $a_{101}=a_{002}=0$ or $c_{200}=c_{110}=c_{020}=0$.
\end{theorem}

\section{Acknowledgments}

The first author is partially supported by S\~ao Paulo Research Foundation (FAPESP) grant 19/10269-3. The second author is supported by S\~ao Paulo Research Foundation (FAPESP) grant 19/13040-7.

\addcontentsline{toc}{chapter}{Bibliografia}
\bibliographystyle{siam}
\bibliography{Referencias.bib}
\end{document}